\documentclass[11pt, a4paper, UKenglish]{article}

\usepackage[UKenglish]{babel}
\usepackage{amssymb, amsmath, textcomp, amsthm}
\usepackage{comment}
\usepackage{bbm,dsfont}
\usepackage{german}
\usepackage[pdftex]{graphicx}
\usepackage[utf8]{inputenc}
\usepackage{authblk}
\numberwithin{equation}{section}
\title{Oscillatory integrals related to Carleson's theorem: fractional monomials}
\author{Shaoming Guo}
\date{}

\def\R{\mathbb{R}}
\def\N{\mathbb{N}}
\def\C{\mathbb{C}}
\def\Z{\mathbb{Z}}

\def\lesim{\lesssim}

\def\beq{\begin{equation}}
\def\endeq{\end{equation}}

\theoremstyle{plain}
\newtheorem{thm}{Theorem}[section]
\newtheorem{prop}[thm]{Proposition}

\newtheorem{lem}[thm]{Lemma}
\newtheorem{cor}[thm]{Corollary}

\newtheorem{rem}[thm]{Remark}
\newtheorem{question}[thm]{Question}

\newtheorem*{conj*}{Conjecture}
\newtheorem*{openproblem*}{Open Problem}

\begin{document}
\maketitle

\begin{abstract}
Stein and Wainger \cite{SW} proved the $L^p$ bounds of the polynomial Carleson operator for all integer-power polynomials without linear term. In the present paper, we partially generalise this result to all fractional monomials in dimension one. Moreover, the connections with Carleson's theorem and the Hilbert transform along vector fields or (variable) curves 
are also discussed in details.
\end{abstract}

\let\thefootnote\relax\footnote{Date: \date{\today}}


\section{Introduction}
In this paper, we will consider two operators related to Carleson's theorem. Fix $\epsilon\in \R$, for a one dimensional Schwartz function $f$, define 
\beq\label{1101ee1.1}
\mathcal{C}_{\epsilon} ^{even}f(x):=\sup_{A\in \R}\left|\int_{\R}e^{i A |y|^{\epsilon}}f(x-y)\frac{dy}{y}\right|.
\endeq
Moreover, define
\beq\label{1101ee1.2}
\mathcal{C}_{\epsilon} ^{odd}f(x):=\sup_{A\in \R}\left|\int_{\R}e^{i A\cdot \text{sgn}(y)\cdot |y|^{\epsilon}}f(x-y)\frac{dy}{y}\right|.
\endeq
The main result we will prove is
\begin{thm}\label{theorem1.1}
For any fixed $\epsilon_1\in \R, \epsilon_1\neq 1$ and $p\in (1, \infty)$, there exists a constant $C_{p, \epsilon_1}>0$ depending on $\epsilon_1$ and $p$ such that
\beq\label{operator}
\left\|\mathcal{C}_{\epsilon_1}^{even}f\right\|_p \le C_{p, \epsilon_1} \|f\|_p.
\endeq
Moreover, for any fixed $\epsilon_2\in \R, \epsilon_2\neq 0$ and $p\in (1, \infty)$, there exists a constant $C_{p, \epsilon_2}>0$ depending on $\epsilon_2$ and $p$ such that
\beq\label{1101ee1.4}
\left\|\mathcal{C}_{\epsilon_2}^{odd}f\right\|_p \le C_{p, \epsilon_2} \|f\|_p.
\endeq
\end{thm}
%
\begin{rem}\label{1101remark1.2}
The estimate \eqref{operator} fails for $\epsilon_1=1$ and the estimate \eqref{1101ee1.4} fails for $\epsilon_2=0$. To see the latter, we just need to notice that 
\beq
\mathcal{C}_0^{odd}f(x)= \left|\int_{\R}f(x-y)\frac{dy}{|y|}\right|,
\endeq
which is clearly not bounded on $L^p$ for any $p\in (1, \infty)$. To see that \eqref{operator} fails for $\epsilon_1=1$, it suffices to show that the corresponding operator without taking supremum, which is 
\beq\label{1101ee1.6}
\int_{\R} e^{i|y|}f(x-y)\frac{dy}{y},
\endeq 
fails to be bounded on $L^2$. By Plancherel's theorem, it is enough to show that the Fourier transform of the convolution kernel in \eqref{1101ee1.6} is unbounded. This calculation can be done explicitly as follows:
\beq
\frac{e^{i|y|}}{y}=\frac{\cos y}{y}+ i \frac{\sin |y|}{y}.
\endeq
The Fourier transform of the real part $\frac{\cos y}{y}$ is a bounded function, hence we just need to consider the imaginary part. 
\beq
\begin{split}
& \mathcal{F}\left(\frac{\sin |y|}{y}\right)=\mathcal{F}\left(\text{sgn}(y)\cdot \frac{\sin y}{y} \right)\\
&=\mathcal{F}(\text{sgn}(y))*\mathcal{F}\left(\frac{\sin y}{y} \right)=\frac{c_0}{\eta}* \chi_{[-1, 1]}(\eta),
\end{split}
\endeq
where $\mathcal{F}$ denotes taking the Fourier transform and $c_0\in \C$ is some numerical constant. It is easy to see that the last term in the above expression is unbounded.\\
\end{rem}

Let us mention some history of the study of the operators \eqref{1101ee1.1} and \eqref{1101ee1.2}. For positive $\epsilon_1$ and $\epsilon_2$, several special cases of the above theorem have already been quite well-known: The case $\epsilon_1=0$ in the estimate \eqref{operator} is the classical Hilbert transform. The case $\epsilon_2=1$ in the estimate \eqref{1101ee1.4} is Carleson's celebrated theorem, and the original proof was given by Carleson \cite{Car}. Later, Fefferman \cite{Fef}, Lacey and Thiele \cite{LT} provided two new proofs.  The cases $\epsilon_1= 2k$ in \eqref{operator} and  $\epsilon_2=2k+1$ in \eqref{1101ee1.4} for  all positive integers $k$ were proven by Stein and Wainger in \cite{SW}. Indeed, the result in \cite{SW} holds true for all integer-power polynomials without linear term. Both the results by Carleson and by Stein and Wainger were unified by Lie in \cite{Lie1} and \cite{Lie2}, where it is proven that for any $d\in \N$, if we denote by $\mathcal{Q}_d$ the class of all integer-power polynomials $Q$ with $\deg(Q)\le d$, then there exists a constant $C_{p, d}>0$ such that
\beq
\left\|\sup_{Q\in \mathcal{Q}_d}\left|\int_{\mathbb{T}}e^{i Q(y)}f(x-y)\frac{dy}{y} \right|\right\|_p\le C_{p, d} \|f\|_p, \forall p\in (1, \infty).
\endeq

\vspace{2mm}

For negative $\epsilon_1$ and $\epsilon_2$, bounds of the form \eqref{operator}  and \eqref{1101ee1.4} are only known for a fixed coefficient $A\in \R$, instead of taking the supremum over all $A\in \R$. Without loss of generality we only consider the estimate \eqref{operator}. By taking $A=1$, we obtain
\beq\label{10Nov1.3}
\int_{\R}e^{i|y|^{\epsilon}}f(x-y)\frac{dy}{y}.
\endeq
The $L^p$ bounds of the above operator for all $p\in (1, \infty)$ were obtained by Hirschman \cite{Hirsch}. The weak type $(1, 1)$ estimate of \eqref{10Nov1.3} is a deep result due to Fefferman \cite{Fefferman}. For the recent development, especially for the generalisation of \eqref{10Nov1.3} from monomial phases to rational phases , see Folch-Gabayet and Wright \cite{FW1}, \cite{FW2} and \cite{FW3}.\\

Our result in Theorem \ref{theorem1.1} should be viewed as a generalisation of the one by Stein and Wainger in \cite{SW}. Indeed, the main tools that we will be using are also essentially the same as those in \cite{SW}, namely the $T T^*$ argument and the stationary phase method. However the techniques that are used by Stein and Wainger in \cite{SW} for the case of $\epsilon$ being an integer do not work for general $\epsilon$. The reason is, that when estimating the kernel of $T T^*$, which is the left hand side of \eqref{EE2.11} below, Stein and Wainger expanded the $\epsilon$-th power polynomial by taking advantage of the fact that $\epsilon$ is an integer, and showed that the phase function always oscillates ``fast'' outside a ``small set''.

%


In our case, i.e. in the case of $\epsilon$ being a general real number (except that $\epsilon=1$), after eliminating Stein and Wainger's ``small set'' (the set $\tilde{\xi}\in (-2^{\theta_1 j}, 2^{\theta_1 j})$ in Page 6 in our case), the phase function in \eqref{EE2.11} might still oscillate very ``slowly''. Hence we need to analyse the phase function more carefully, which is done in our crucial Lemma \ref{lemma2.2}.

The novelty of Lemma \ref{lemma2.2} is, that except for eliminating the ``small set'' by Stein and Wainger, we need to eliminate another ``small set'', which is the set in \eqref{EE3.23}, and only by doing this will the phase function in \eqref{EE2.11} oscillate ``fast''.

\subsection{Uniform estimates}
A slight modification of the proof of Theorem \ref{theorem1.1} leads to the following uniform estimates:
\begin{thm}\label{1401theorem1.3}
For any $0<\delta<1$, there exists a constant $C_{\delta}$ such that for all $\epsilon$ with $|\epsilon|>\delta$ and $|\epsilon-1|>\delta$, we have 
\beq\label{1101ee1.11}
\left\|\mathcal{C}_{\epsilon}^{even}f\right\|_2 \le C_{\delta} \|f\|_2,
\endeq
and 
\beq\label{1101ee1.12}
\left\|\mathcal{C}_{\epsilon}^{odd}f\right\|_2 \le C_{\delta} \|f\|_2.
\endeq
\end{thm}

The argument in Remark \ref{1101remark1.2} indicates that the estimate \eqref{1101ee1.11} blows up when $\epsilon$ tends to 1, and the estimate \eqref{1101ee1.12} blows up when $\epsilon$ tends to 0. However, there are still two other blow-ups, namely when $\epsilon\to 0$ in \eqref{1101ee1.11} and $\epsilon\to 1$ in \eqref{1101ee1.12} separately. For the former case, letting $\epsilon\to 0$, we obtain a ``limit''
\beq
C_0^{even}f(x)=\left|\int_{\R}f(x-y)\frac{dy}{y} \right|,
\endeq
which is the classical Hilbert transform. For the latter case, letting $\epsilon\to 1$ in \eqref{1101ee1.12}, we obtain a ``limit''
\beq\label{1201ee1.14}
C_1^{odd}f(x)=\sup_{A\in \R}\left|\int_{\R}f(x-y)e^{iAy}\frac{dy}{y} \right|,
\endeq
which is exactly Carleson's maximal operator. Hence it is reasonable to ask the following
\begin{question}
Is there a universal constant $C_0>0$ such that for all $|\epsilon_1|<1/2$ and all $|\epsilon_2-1|<1/2$, we have 
\beq\label{1101ee1.15}
\left\|\mathcal{C}_{\epsilon_1}^{even}f\right\|_2 \le C_{0} \|f\|_2,
\endeq
and 
\beq\label{1101ee1.16}
\left\|\mathcal{C}_{\epsilon_2}^{odd}f\right\|_2 \le C_{0} \|f\|_2?
\endeq
\end{question}

\begin{rem}
If the estiamte \eqref{1101ee1.16} were true, then by a simple limiting argument, it would imply Carleson's theorem.
\end{rem}

To support the above question, we prove the following uniform $L^{\infty}$ estimate for the multiplier of the convolution kernels $e^{i|y|^{\epsilon}}/y$ and  $e^{i\cdot \text{sgn}(y)\cdot |y|^{\epsilon}}/y$. This will imply the uniform $L^2$ boundedness of the operators \eqref{1101ee1.1} and \eqref{1101ee1.2} without taking the supremum over $A\in \R$.

\begin{thm}\label{1401theorem1.6}
There exists a universal constant $C_0>0$ such that for all $|\epsilon_1|<1/2$, we have 
\begin{equation}\label{1401eee1.17}
\left\| \int_{\R}e^{i|t|^{\epsilon_1}} e^{-i\lambda t}\frac{dt}{t}\right\|_{L^{\infty}(\lambda)} \le C_{0},
\end{equation}
and for all $1/2<\epsilon_2<3/2$, we have 
\begin{equation}\label{1401eee1.18}
\left\| \int_{\R}e^{i\cdot \text{sgn}(t)\cdot |t|^{\epsilon_2}} e^{-i\lambda t}\frac{dt}{t}\right\|_{L^{\infty}(\lambda)} \le C_{0}.
\end{equation}
\end{thm}

\subsection{Connection with the Hilbert transform along vector fields or (variable) curves}

The results in Theorem \ref{theorem1.1} are closely related to the Hilbert transform along planar vector fields or curves. We start with the case of planar vector fields. Before explaining the relation, we state the following result due to Bateman and Thiele \cite{BT} concerning the $L^p$ bounds of the Hilbert transform along the so-called one-variable vector fields.

\begin{thm}(\cite{Bateman}, \cite{BT})\label{1201theorem1.7}
Let $u:\R\to \R$ be an arbitrary measurable function. Then the Hilbert transform along the one-variable vector fields $(1, u): \R^2\to \R^2$, which is given by
\beq\label{1201ee1.19}
H^u_1 f(x_1, x_2):=\int_{\R}f(x_1-t, x_2-u(x_1)t)\frac{dt}{t},
\endeq
is bounded on $L^p$ for all $p>3/2$.
\end{thm}

The case $p=2$ in the above Theorem \ref{1201theorem1.7} is very special as it is equivalent with the $L^2$ bounds of Carleson's maximal operator \eqref{1201ee1.14}. This was first observed by Coifman and El Kohen. We review the discussion as presented in \cite{BT}. Denoting by $\hat{f}$ the partial Fourier transform in the $x_2$ variable on the plane we obtain formally
\begin{equation}\label{onevarvf}
 \int f(x_1-t, x_2-u(x_1)t)\frac{dt}{t}
\end{equation}
$$ = \int e^{i x_2 \xi_2} \int \widehat{f}(x_1-t, \xi_2)
e^{i u(x_1)t\xi_2} \frac{dt}{t}d\xi_2.$$
By the Plancherel theorem, 
\beq\label{NN1.8}
\|H_{1}^u f\|_2=\|\int \widehat{f}(x_1-t, \xi_2)
e^{i u(x_1)t\xi_2} \frac{dt}{t}\|_2
\endeq
For each fixed $\xi_2$, we recognize this to
essentially be the linearisation of Carleson's maximal operator \eqref{1201ee1.14}. Hence the right hand side of \eqref{NN1.8} can be bounded by
\beq\label{NN1.10}
\|C_{1}^{odd}\hat{f}(x_1, \xi_2)\|_2\le C \|\hat{f}(x_1, \xi_2)\|_2 \le C \|f\|_2,
\endeq
for some positive constant $C\in \R$. Moreover, by choosing the function $u$ properly in \eqref{1201ee1.19}, the $L^2$ boundedness of $H_1^u$ also implies the $L^2$ boundedness of Carleson's maximal operator. \\

In the same way that the $L^2$ bounds of Carleson's maximal operator imply the $L^2$ bounds of $H_1^{u}$ in \eqref{1201ee1.19}, the results in Theorem \ref{theorem1.1} has the following corollary concerning the $L^2$ bounds of the Hilbert transform along certain variable curves. 

\begin{cor}\label{1101cor1.6}
Let $u:\R\to \R$ be an arbitrary measurable function. Fix $\epsilon\in \R$, define
\beq
H_{\epsilon}^u f(x_1, x_2):=\int_{\R}f(x_1-t, x_2-u(x_1)\text{sgn}(t)\cdot |t|^{\epsilon})\frac{dt}{t}.
\endeq
Then for any $\epsilon\neq 0$, there exists $C_{\epsilon}>0$ such that 
\beq\label{1401eee1.24}
\|H_{\epsilon}^u f\|_2 \le C_{\epsilon} \|f\|_2.
\endeq
\end{cor}
\begin{rem}
A similar result holds true for the Hilbert transform along the even curve $(t, |t|^{\epsilon})$.
\end{rem}

However, the above argument by Coifman and El Kohen works only in $L^2$. So far it is not know whether the result in Corollary \ref{1101cor1.6} can be generalised to any $p$ other than 2.\\


The result in Corollary \ref{1101cor1.6} is a generalisation of the $L^2$ boundedness of the Hilbert transform along a fixed odd curve $(t, \text{sgn}(t)\cdot |t|^{\epsilon})$ (or an even curve $(t, |t|^{\epsilon})$). The Hilbert transform along curve $(t, \gamma(t))$ for some $\gamma:\R \to \R$, which is defined as
\beq
H_{\gamma}f(x_1, x_2)=\int_{\R}f(x_1-t, x_2-\gamma(t))\frac{dt}{t},
\endeq
has been extensively studied, see for example \cite{CCC}, \cite{NRW1}, \cite{NRW2}, \cite{NVWW1} and \cite{NVWW2}. Here we only state the following
\begin{thm}(\cite{CCC})\label{1401theorem1.10}
Let $\gamma:\R\to \R$ be an even and convex function with $\gamma(0)=\gamma'(0)=0$. Then for all $p\in (1, \infty)$, a necessary and sufficient condition for the boundedness of $H_{\gamma}$ on $L^p$ is that 
\beq\label{1401eee1.26}
\text{there exists } \lambda>1 \text{ with } \gamma'(\lambda t)\ge 2 \gamma'(t) \text{ for all } t>0.
\endeq
\end{thm}

By comparing the result in Theorem \ref{1401theorem1.10}  with the one in Corollary \ref{1101cor1.6}, it might be reasonable to expect that the estimate \eqref{1401eee1.24} holds true for a larger class of curves satisfying conditions like \eqref{1401eee1.26}.\\

{\bf Organisation of Paper:} In Section \ref{1401section2} we will present the proof of Theorem \ref{theorem1.1}. The main argument is based on the $T T^{*}$ method and the oscillatory integral estimates. 

In Section \ref{1602section3} we will prove the uniform estimates in Theorem \ref{1401theorem1.3}. The proof is a slight modification of the one for Theorem \ref{theorem1.1}. 

In the last Section \ref{1602section4} we will prove the uniform estimates in Theorem \ref{1401theorem1.6} concerning the $L^{\infty}$ bounds of the Fourier transform of certain convolution kernels. The proof is based on careful integration by parts. \\

{\bf Notations:} Throughout this paper, we will write $x\ll y$ to mean that $x\le y/10$, $x\lesim y$ to mean that there exists a universal constant $C$ s.t. $x\le C y$, and $x\sim y$ to mean that $x\lesim y$ and $y\lesim x$.  $\mathbbm{1}_E$ will always denote the characteristic function of the set $E$.\\

%
%
%
%
%
%
%
%
%
%
%
%
%
%
%

{\bf Acknowledgements.} The author would like to thank his advisor, Prof. Christoph Thiele, for helpful discussions. The author also thanks Prof. Po Lam Yung for his valuable comments.

\section{Proof of Theorem \ref{theorem1.1}}\label{1401section2}

It turns out that in the following proof of Theorem \ref{theorem1.1}, there are no distinguished differences between the cases $\epsilon_i>0$ and $\epsilon_i<0$ (here $i=1, 2$). Therefore in most part of this section we will be talking about the case $\epsilon_i>0$ (which is also slightly more tricky as it includes the threshold $\epsilon_i=1$), and leave the discussion of the case $\epsilon_i<0$ till the end as a remark. \\

The structure of this section is as follows. In the first subsection we state the strategy of the proof of the case $\epsilon_i>0$ in Theorem \ref{theorem1.1}. The main idea is that we first decompose the operator on the left hand side of \eqref{operator} or \eqref{1101ee1.4} into two parts: the high frequency part and the low frequency part (see the following \eqref{10nov2.4}). The high frequency part will be dominated pointwise by the maximal operator and the maximal Hilbert transform (see Lemma \ref{lowfrequency}). For the $L^2$ bounds of the low frequency part, we will apply the $TT^{*}$ method and techniques from oscillatory integrals to obtain certain exponential decay (see Proposition \ref{prop2.2}). The $L^p$ bounds of the low frequency part follow simply by interpolating the $L^2$ bounds with certain trivial bounds. In the second subsection, we will give the details of the proof of Proposition \ref{prop2.2}. In the last subsection, we will remark on the proof of the case $\epsilon_i<0$ in Theorem \ref{theorem1.1}.

\subsection{Strategy of the proof of Theorem \ref{theorem1.1} for $\epsilon_i>0$}\label{section2}

The proofs for \eqref{operator} and \eqref{1101ee1.4} in Theorem \eqref{theorem1.1} are similar, hence here we will only consider the former case. After a linearisation of the maximal operator on the left hand side of \eqref{operator}, we are going to prove the boundedness of 
\beq\label{10nov2.1}
T_A f(x):=\int_{\R}e^{i A(x) |y|^{\epsilon}}f(x-y)\frac{dy}{y},
\endeq
with a bound being independent of the positive measurable function $A:\R\to \R^+.$ \\

Take a smooth partition of unity
\beq\label{1401ee2.2}
\sum_{j=-\infty}^{\infty}\psi_j(y)=1, \forall y\neq 0,
\endeq
where $\psi_j(y)=\psi(2^{j}y)$. Then it is easy to see that
\beq
\sum_{j=-\infty}^{\infty}\psi_j(A(x)^{1/\epsilon}y)=1, \forall x\in \R.
\endeq
Hence we can split our operator $T_A$ into the following two parts:
\beq\label{10nov2.4}
T_A f(x)=\left(\sum_{j> 0}+\sum_{j\le 0}\right) \int_{\R} e^{i A(x) |y|^{\epsilon}} \psi_j(A(x)^{1/\epsilon}y)\frac{f(x-y)}{y}dy.
\endeq
The same decomposition has already been used in \cite{DL}. For the former part, we denote  it as
\beq\label{10nov2.5}
T_A^{high}f(x):=\sum_{j> 0} \int_{\R} e^{i A(x) |y|^{\epsilon}} \psi_j(A(x)^{1/\epsilon}y)\frac{f(x-y)}{y}dy,
\endeq
while for the latter part, we denote it as
\beq\label{GG2.6}
T_A^{low}f(x):=\sum_{j\le 0} \int_{\R} e^{i A(x) |y|^{\epsilon}} \psi_j(A(x)^{1/\epsilon}y)\frac{f(x-y)}{y}dy.
\endeq

The boundedness of the former part is done in the following 
\begin{lem}\label{lowfrequency}
Under the above notations, we have the following pointwise estimate
\beq
|T_A^{high}f(x)|\lesim M f(x)+ H^* f(x),
\endeq
where $M$ denotes the one-dimensional Hardy-Littlewood maximal operator, and $H^*$ denotes the maximal Hilbert transform.
\end{lem}
\noindent {\bf Proof of Lemma \ref{lowfrequency}:} The idea is to approximate the term $e^{iA(x)|y|^{\epsilon}}$ by $1$, as in the operator $T_A^{high}$, the exponent $A(x)|y|^{\epsilon}$ is always small. By denoting
\beq
\phi_0:=\sum_{j> 0}\psi_j,
\endeq
we obtain
\beq
T_A^{high}f(x)=\int_{\R} e^{i A(x) |y|^{\epsilon}} \phi_0(A(x)^{1/\epsilon}y)\frac{f(x-y)}{y}dy.
\endeq
By subtracting a zero we obtain
\beq
\begin{split}
T_A^{high}f(x)&=\int_{\R} \left(e^{i A(x) |y|^{\epsilon}}-1\right) \phi_0(A(x)^{1/\epsilon}y)\frac{f(x-y)}{y}dy\\
		&+\int_{\R} \phi_0(A(x)^{1/\epsilon}y)\frac{f(x-y)}{y}dy.
\end{split}
\endeq
For the latter part, we bound it by the maximal Hilbert transform, i.e. 
\beq
\left|\int_{\R} \phi_0(A(x)^{1/\epsilon}y)\frac{f(x-y)}{y}dy\right|\lesim H^*f(x).
\endeq
For the former part, we will bound it by the maximal operator, i.e. 
\beq
\begin{split}
& \left|\int_{\R} \left(e^{i A(x) |y|^{\epsilon}}-1\right) \phi_0(A(x)^{1/\epsilon}y)\frac{f(x-y)}{y}dy\right|\\
& \lesim \int_{\R} A(x) |y|^{\epsilon} \phi_0(A(x)^{1/\epsilon}y)\frac{|f(x-y)|}{|y|}dy \lesim Mf(x).
\end{split}
\endeq
So far we have finished the proof of Lemma \ref{lowfrequency}.$\Box$\\

Hence what is left is to prove the $L^p$ boundedness of the low frequency part, i.e. the expression in \eqref{GG2.6}. If we denote
\beq
T^j f(x):=\int_{\R} e^{i A(x) |y|^{\epsilon}} \psi_j(A(x)^{1/\epsilon}y)\frac{f(x-y)}{y}dy,
\endeq
then
\beq
T_A^{low}f(x)=\sum_{j\le 0}T^j f(x).
\endeq
Hence by the triangle inequality, it suffices to prove the following
\begin{prop}\label{prop2.2}
Fix $\epsilon\neq 1$ and $p\in (1, \infty)$. For any non-positive integer $j$, for an arbitrary positive measurable function $A$, we have
\beq\label{NN2.15}
\|T^j f\|_p \lesim 2^{\sigma j}\|f\|_p,
\endeq
with $\sigma>0$ being independent of $A$ and $j$.
\end{prop}

To prove Proposition \ref{prop2.2},  we first observe that 
\beq
|T^j f(x)|\le \int_{\R} |\psi_j(A(x)^{1/\epsilon}y)\frac{f(x-y)}{y}|dy \lesim M f(x),
\endeq
which then implies the trivial bound
\beq
\|T^j f\|_p \lesim \|M f\|_p \lesim \|f\|_p.
\endeq
Hence to prove \eqref{NN2.15} for an arbitrary $p>1$, by interpolation, we just need to prove the case $p=2$.

%
%
%
%
%
%
%
%
%
%
%
%
%
%
%
%
%
%
%
%
%

\subsection{Proof of the Proposition \ref{prop2.2} for $p=2$}\label{section3}
\subsubsection{Calculating the kernel of $T T^*$}
%

To obtain the $L^2$ bounds in \eqref{NN2.15}, we want to use the $T T^*$ method. First, we write down the dual operator, which is  
\beq\label{1401eee2.18}
T^{j,*}g(y)=\int_{\R} e^{-i A(x)|x-y|^{\epsilon}}\psi_j\left(A(x)^{1/\epsilon}(x-y)\right)\frac{g(x)dx}{x-y}.
\endeq
Therefore,
\beq
\begin{split}
& T^j T^{j,*}f(y) = T^j\left(\int_{\R} e^{-i A(x)|x-y|^{\epsilon}} \psi_j\left(A(x)^{1/\epsilon}(x-y)\right)\frac{ f(x)}{x-y}dx\right)\\
		&= \int_{\R}\int_{\R}e^{i A(y)|y-z|^{\epsilon}} \frac{\psi_j\left(A(y)^{1/\epsilon}(y-z)\right)}{y-z} e^{-i A(x)|x-z|^{\epsilon}}\frac{\psi_j\left(A(x)^{1/\epsilon}(x-z)\right)}{x-z}dz f(x)dx\\
		&= \int_{\R}\int_{\R}e^{iA(y)|y-x+z|^{\epsilon}} \frac{\psi_j\left(A(y)^{1/\epsilon}(y-x+z)\right)}{y-x+z} e^{-i A(x)|z|^{\epsilon}}\frac{\psi_j\left(A(x)^{1/\epsilon}z\right)}{z}dz f(x)dx\\
		&=: \int_{\R}(\Phi_j^{A(y)}*\tilde{\Phi}_j^{A(x)})(y-x)f(x)dx,
\end{split}
\endeq
where
\beq
\Phi_j^{A}(\xi):=e^{iA |\xi|^{\epsilon}}\frac{\psi_j(A^{1/\epsilon}\xi)}{\xi},
\endeq
and 
\beq
\tilde{\Phi}_j^A (\xi):=\bar{\Phi}_j^A(-\xi).
\endeq

\vspace{4mm}


In the following calculation, we assume w.l.o.g. that $A(x)\le A(y)$. If we denote $\xi=-y+x$, then the kernel of the operator $T^j T^{j, *}$ is given by
\beq\label{NN3.5}
\begin{split}
& \Phi_j^{A(y)}*\tilde{\Phi}_j^{A(x)}(\xi)\\
&= \int_{\R}e^{i A(y)|\eta|^{\epsilon}}\frac{\psi_0(2^j A(y)^{\frac{1}{\epsilon}}\eta)}{\eta} e^{-i A(x)|\xi-\eta|^{\epsilon}}\frac{\psi_0(2^j A(x)^{\frac{1}{\epsilon}}(\xi-\eta))}{\xi-\eta} d\eta.
\end{split}
\endeq
To evaluate the above integral, we do the following change of variable
\beq
2^j A(y)^{\frac{1}{\epsilon}}\eta=\eta,
\endeq
and denote
\beq
\left(\frac{A(x)}{A(y)}\right)^{\frac{1}{\epsilon}}=h,
\endeq
%
then the expression in \eqref{NN3.5} becomes
\beq
2^j A(x)^{\frac{1}{\epsilon}} \int_{\R}e^{i |\eta|^{\epsilon}2^{-j\epsilon}} \frac{\psi_0(\eta)}{\eta} e^{-i |A(x)^{\frac{1}{\epsilon}}\xi-2^{-j} h \eta|^{\epsilon}}\frac{\psi_0(2^j A(x)^{\frac{1}{\epsilon}}\xi- h \eta)}{2^j A(x)^{\frac{1}{\epsilon}}\xi- h \eta}  d\eta.
\endeq
If we further denote
\beq
2^j A(x)^{\frac{1}{\epsilon}}\xi=:\tilde{\xi},
\endeq
we then obtain
%
%
\beq\label{EE3.15}
\begin{split}
2^j A(x)^{\frac{1}{\epsilon}}  \int_{\R}e^{i |2^{-j} \eta|^{\epsilon}-i 2^{-j\epsilon}|\tilde{\xi}-h \eta|^{\epsilon}} \frac{\psi_0(\eta)}{\eta} \frac{\psi_0(\tilde{\xi}- h \eta)}{\tilde{\xi}- h \eta} d\eta.
\end{split}
\endeq

As the next step, we will prove
\begin{lem}\label{lemma2.1}
There exists two small positive real numbers $\lambda_1$ and $\lambda_2$ such that the following pointwise estimate in $\tilde{\xi}$ holds
\beq\label{EE2.11}
\begin{split}
& \left|\int_{\R}e^{i |2^{-j} \eta|^{\epsilon}-i 2^{-j\epsilon}|\tilde{\xi}-h \eta|^{\epsilon}} \frac{\psi_0(\eta)}{\eta} \frac{\psi_0(\tilde{\xi}- h \eta)}{\tilde{\xi}- h \eta} d\eta \right| \\
& \lesim \chi_{[-2^{\lambda_1 j}, 2^{\lambda_1 j}]}(\tilde{\xi})+ 2^{\lambda_2 j}\chi_{[-3, 3]}(\tilde{\xi}),
\end{split}
\endeq
where the constant is independent of $h\in (0, 1]$.
\end{lem}

Hence by Stein and Wainger's small set maximal function theorem (Proposition 3.1 in \cite{SW}), we obtain that there exists $\sigma(\lambda_1, \lambda_2)>0$ which depends on $\lambda_1$ and $\lambda_2$ such that
\beq
\|T^j f\|_2 \lesim 2^{\sigma(\lambda_1, \lambda_2)j}\|f\|_2.
\endeq

To prove the above lemma, we need to analyse the phase function on the left hand side of \eqref{EE2.11} carefully. As the functions $|\eta|^{\epsilon}$ and $|\tilde{\xi}-h\eta|^{\epsilon}$ are not smooth due to the fact that we are taking the absolute values of $\eta$ and $\tilde{\xi}-h\eta$, we need to divide the analysis into four cases. First w.l.o.g. we assume that $\eta>0$, i.e. we are taking the positive branch of $\psi_0$.  Then we denote by {\bf Case One} the case when
\beq
\tilde{\xi}-h\eta>0
\endeq
in the term $|\tilde{\xi}-h\eta|^{\epsilon}$ on the left hand side of \eqref{EE2.11}, and {\bf Case Two} the case when
\beq
\tilde{\xi}-h\eta<0.
\endeq

\subsubsection{Proof of Lemma \ref{lemma2.1}: Case One}

This case is the easier case, as we will see that the phase function
\beq
\Phi_h(\tilde{\xi},\eta):=2^{-\epsilon j}\left(\eta^{\epsilon}-  \left(\tilde{\xi} - h\eta\right)^{\epsilon}\right)
\endeq
in \eqref{EE2.11} will always oscillate fast, which means we can apply the stationary phase method directly. \\

Under the assumption that 
\beq
\tilde{\xi}-h\eta>0,
\endeq
the derivative of the phase function $\Phi_h(\tilde{\xi},\eta)$ becomes
\beq
\epsilon \cdot  2^{-\epsilon j}\left(\eta^{\epsilon-1}+ h \left(\tilde{\xi} - h\eta\right)^{\epsilon-1}\right) \gtrsim 2^{-\epsilon j},
\endeq
which means the derivative of the phase function is bounded from below. \\

However, to apply the stationary phase method, we still need the derivative of the phase function to be monotone, which is not always the case. Fortunately, the second order derivative of the phase function has only one critical point, hence by Proposition 2 in Page 332 of Stein's book \cite{Stein}, we obtain
\beq\label{EE2.18}
\left|\int_{\R}e^{i (2^{-j} \eta)^{\epsilon}-i 2^{-j\epsilon}(\tilde{\xi}-h \eta)^{\epsilon}} \frac{\psi_0(\eta)}{\eta} \frac{\psi_0(\tilde{\xi}- h \eta)}{\tilde{\xi}- h \eta} d\eta \right| \lesim 2^{c_1 j}\chi_{[-3, 3]}(\tilde{\xi}),
\endeq
for some positive real $c_1$, which finishes the proof of Case One.

\subsubsection{Proof of Lemma \ref{lemma2.1}: Case Two}

This time the derivative of the phase function $\Phi_h(\tilde{\xi},\eta)$ is given by
\beq\label{GG3.19}
\epsilon \cdot 2^{-\epsilon j}\left(\eta^{\epsilon-1}- h\cdot  \left(h\eta-\tilde{\xi}\right)^{\epsilon-1}\right).
\endeq
The analysis of this term is a bit more involved than the one in the last case, as in \eqref{GG3.19} there is not only turning points of the derivative of the phase function, but also turning points of the phase function, which means for certain choices of $h$ and $\tilde{\xi}$, we are not allowed to use the stationary phase method.\\

{\bf Case $0< h\le h_0:$} here $h_0$ is some small positive number to be chosen later. As $h$ is small, we observe that
\beq
\eta^{\epsilon-1}- h\cdot  (h\eta-\tilde{\xi})^{\epsilon-1}
\endeq
from \eqref{GG3.19} is monotone in $\eta$ and is bounded from below by some constant, say $1/10$. Hence the stationary phase method applies and we again obtain an estimate of the form \eqref{EE2.18}.\\


{\bf Case $h_0< h\le 1$ and $\tilde{\xi}\in (-2^{\theta_1 j}, 2^{\theta_1 j}):$} here $\theta_1$ is some small positive number to be determined. In this case, we use the trivial bound, i.e. to bound the exponential factor
\beq
\exp[i \Phi_h(\tilde{\xi}, \eta)]=\exp[i (2^{-j} \eta)^{\epsilon}-i 2^{-j\epsilon}(\tilde{\xi}-h \eta)^{\epsilon}]
\endeq 
by one. Hence 
\beq
\left|\int_{\R}e^{i (2^{-j} \eta)^{\epsilon}-i 2^{-j\epsilon}(\tilde{\xi}-h \eta)^{\epsilon}} \frac{\psi_0(\eta)}{\eta} \frac{\psi_0(\tilde{\xi}- h \eta)}{\tilde{\xi}- h \eta} d\eta \right| \lesim \chi_{[-2^{\theta_1 j}, 2^{\theta_1 j}]}(\tilde{\xi}),
\endeq
which is another term on the right hand side of \eqref{EE2.11} with $\lambda_1:=\theta_1$.\\

{\bf case $h_0< h\le 1$ and $|\tilde{\xi}|\ge 2^{\theta_1 j}:$} the derivative of the phase function might be small in this case, hence we single out the set
\beq\label{EE3.23}
E_{h, \tilde{\xi}}:=\{\eta: \eta\in (1/2, 5/2), h\eta-\tilde{\xi}\in (1/2, 5/2) \text{ and } |\eta^{\epsilon-1}-h\cdot (h \eta-\tilde{\xi})^{\epsilon-1}|\le 2^{\theta_2 j} \},
\endeq
where $\theta_2$ is a positive number smaller than $\epsilon$. Hence outside the set $E_{h, \tilde{\xi}}$, the derivative of the phase function has a lower bound, i.e.
\beq
\eqref{GG3.19} \ge 2^{-(\epsilon-\theta_2)j}.
\endeq
Moreover we have the following crucial upper bound on the size of the bad set
%
\begin{lem}\label{lemma2.2}
There exists a constant $c_2>0$ which depends on $\theta_1$ and $\theta_2$ such that
\beq
|E_{h, \tilde{\xi}}|\lesim 2^{c_2 j}.
\endeq
\end{lem}

We postpone the proof of Lemma \ref{lemma2.2} to the next subsection and proceed with the proof of Lemma \ref{lemma2.1}. On the bad set $E_{h, \tilde{\xi}}$, we apply the estimate in Lemma \ref{lemma2.2} and the same trivial bound as in the last case to obtain
\beq
\left|\int_{E_{h,\tilde{\xi}}}e^{i (2^{-j} \eta)^{\epsilon}-i 2^{-j\epsilon}(\tilde{\xi}-h \eta)^{\epsilon}} \frac{\psi_0(\eta)}{\eta} \frac{\psi_0(\tilde{\xi}- h \eta)}{\tilde{\xi}- h \eta} d\eta \right| \lesim 2^{c_2 j}\chi_{[-3, 3]}(\tilde{\xi}).
\endeq

Outside the bad set $E_{h, \tilde{\xi}}$ we apply the stationary phase principle. Notice that the derivative of the phase function is not monotone, but again there exists only finitely many turning points. Hence we can still apply Proposition 2 in Page 332 in Stein's book \cite{Stein} to obtain the bound
\beq
\left|\int_{\R\setminus E_{h,\tilde{\xi}}}e^{i (2^{-j} \eta)^{\epsilon}-i 2^{-j\epsilon}(\tilde{\xi}-h \eta)^{\epsilon}} \frac{\psi_0(\eta)}{\eta} \frac{\psi_0(\tilde{\xi}- h \eta)}{\tilde{\xi}- h \eta} d\eta \right| \lesim 2^{c_3 j}\chi_{[-3, 3]}(\tilde{\xi}),
\endeq
for some constant $c_3$ depending on $\epsilon, \theta_1$ and $\theta_2$, where $\theta_1$ will be chosen accordingly in the proof of Lemma \ref{lemma2.2}. 

In the end, we just need to take $\lambda_2:=\min\{c_1, c_2, c_3\}$.

\subsubsection{Proof of the crucial Lemma \ref{lemma2.2}} 
The proof is basic, and the main point here is how to make full use of the condition that $|\tilde{\xi}|\ge 2^{\theta_1 j}$. 

As we assume that 
\beq
\eta\in (1/2, 5/2), h\eta-\tilde{\xi}\in (1/2, 5/2) \text{ and } h_0<h\le 1,
\endeq 
by the mean value theorem, we obtain that 
\beq
\left|\eta^{\epsilon-1}-h\cdot (h \eta-\tilde{\xi})^{\epsilon-1}\right| \sim \left| \eta-h^{\frac{\epsilon}{\epsilon-1}}\eta+ h^{\frac{1}{\epsilon-1}}\tilde{\xi}\right|,
\endeq
where the constant depends only on $h_0$ and $\epsilon$. Hence the restriction 
\beq
\left|\eta^{\epsilon-1}-h\cdot (h \eta-\tilde{\xi})^{\epsilon-1}\right| \le 2^{\theta_2 j}
\endeq
in the definition of the set $E_{h,\tilde{\xi}}$ turns to 
\beq\label{EE2.30}
\left| \eta-h^{\frac{\epsilon}{\epsilon-1}}\eta+ h^{\frac{1}{\epsilon-1}}\tilde{\xi}\right|\lesim 2^{\theta_2 j},
\endeq
which further implies that
\beq\label{EE2.31}
|E_{h,\tilde{\xi}}|\lesim 2^{\theta_2 j}(1-h^{\frac{\epsilon}{\epsilon-1}})^{-1}.
\endeq

To control the right hand side of the last expression, the idea is to show that $h$ can not be too close to $1$, due to the restriction that $\tilde{\xi}\ge 2^{\theta_1 j}$. \\

{\bf Case $\epsilon>1$:} we will choose $\theta_1$ such that
\beq\label{EE2.32}
\theta_1=\theta_2/4,
\endeq
and then show that in order for the set $E_{h,\tilde{\xi}}$ not to be empty, we must have
\beq\label{EE2.33}
h<1-2^{\frac{\theta_2}{2}\cdot j}.
\endeq
We argue by contradiction: assume that
\beq
1-2^{\frac{\theta_2}{2}\cdot j}\le h\le 1,
\endeq
then 
\beq
\left| \eta-h^{\frac{\epsilon}{\epsilon-1}}\eta\right|\lesim 2^{\frac{\theta_2}{2}\cdot j}.
\endeq
Hence by the choice of $\theta_1$ and $\theta_2$ in \eqref{EE2.32}, we obtain
\beq
\left| \eta-h^{1+\frac{1}{\epsilon}}\eta+ h^{\frac{1}{\epsilon}}\tilde{\xi}\right| \ge \left|h^{\frac{1}{\epsilon}}\tilde{\xi}\right|-\left| \eta-h^{1+\frac{1}{\epsilon}}\eta\right|\gtrsim 2^{\theta_1 j},
\endeq
which is a contradiction to \eqref{EE2.30}. 

Thus we have verified \eqref{EE2.33}. By substituting \eqref{EE2.33} into the right hand side of \eqref{EE2.31}, we obtain that
\beq
|E_{h,\tilde{\xi}}|\lesim 2^{\frac{\theta_2}{2}\cdot j}.
\endeq

\vspace{4mm}

{\bf Case $\epsilon<1$:} This case is similar to the previous one. We just need to notice that $\frac{\epsilon}{\epsilon-1}<0$, hence instead of \eqref{EE2.33}, what we need to show is
\beq
\frac{1}{h}> 1+2^{\frac{\theta_2}{2}\cdot j}
\endeq
in order for the set $E_{h, \tilde{\xi}}$ not to be empty. The proof is again by a similar contradiction argument as before, hence we leave it out. Thus we have finished the proof of Lemma \ref{lemma2.2}.

\subsection{Remarks on the case $\epsilon<0$}\label{section4}

As has been mentioned before, the proof of Theorem \ref{theorem1.1} for the case $\epsilon<0$ is essentially the same as that for the case $\epsilon>0$, with just minor modifications that we will state in this subsection.\\

Consider the linearised operator \eqref{10nov2.1} for some $\epsilon<0$, in order to distinguish from the case $\epsilon>0$, we replace $\epsilon$ by $-\epsilon$ and denote
\beq
T_A f(x):=\int_{\R}e^{i \frac{A(x)}{|y|^{\epsilon}} }f(x-y)\frac{dy}{y}.
\endeq

The starting point is the same as before, which is to do the high-low frequency decomposition \eqref{10nov2.4} and write
\beq
T_A f(x)=\left(\sum_{j> 0}+\sum_{j\le 0}\right) \int_{\R} e^{i \frac{A(x)}{ |y|^{\epsilon}}} \psi_j\left(\frac{y}{A(x)^{1/\epsilon}}\right)\frac{f(x-y)}{y}dy.
\endeq
Now, instead of \eqref{10nov2.5} and \eqref{GG2.6}, we denote
\beq
T_A^{high} f(x)=\sum_{j\le 0} \int_{\R} e^{i \frac{A(x)}{ |y|^{\epsilon}}} \psi_j\left(\frac{y}{A(x)^{1/\epsilon}}\right)\frac{f(x-y)}{y}dy,
\endeq
and 
\beq
T_A^{low} f(x)=\sum_{j> 0} \int_{\R} e^{i \frac{A(x)}{ |y|^{\epsilon}}} \psi_j\left(\frac{y}{A(x)^{1/\epsilon}}\right)\frac{f(x-y)}{y}dy.
\endeq
Then Lemma \ref{lowfrequency} and Proposition \ref{prop2.2} will stay true by similar arguments. We leave out the details.

\section{Proof of the uniform estimate in Theorem \ref{1401theorem1.3}}\label{1602section3}

Again we will only consider the estimate \eqref{1101ee1.11}, as the proof for the other estimate is similar. The exponent $\epsilon$ lies in the region
\beq
(-\infty, \-\delta]\cup [\delta, 1-\delta]\cup [1+\delta, \infty).
\endeq
The middle part $[\delta, 1-\delta]$ is a closed interval, and the argument in Section \ref{1401section2} can be easily checked to be uniform for $\epsilon$ on this interval. Hence we will need to prove a uniform estimate for $\epsilon$ on the union of intervals 
\beq
(-\infty, \-\delta]\cup [1+\delta, \infty).
\endeq
Here we will carry out the calculation for the case $\epsilon\ge 1+\delta$. The argument for the case $\epsilon\le -\delta$ is similar.\\

We start with the proof. As we are to prove a uniform estimate for large $\epsilon$, we will denote $n:=\epsilon$ to indicate that $\epsilon$ is a large number. Similar to the linearisation done in \eqref{10nov2.1} of Subsection \ref{section2}, it suffices to consider 
\beq
T_A f(x):=\int_{\R} e^{iA(x)|t|^n}f(x-t)\frac{dt}{t},
\endeq
where $A:\R\to \R^+$. The proof below is a slight modification of the one for Theorem \ref{theorem1.1} in Section \ref{1401section2}. However, we need to be careful with the scale of the dyadic decomposition that we do in \eqref{1401ee2.2} as otherwise the bound will blow up when $n\to \infty$.\\

Denote $\lambda:=2^{1/n}$. Choose a smooth function $\psi_0$ which is supported on $(1/\lambda, \lambda^2)\cup (-\lambda^2, -1/\lambda)$ such that 
\beq
\psi_0(t)=1, \forall t\in [1, \lambda]\cup [-\lambda, -1],
\endeq
and
\beq
\sum_{j\in \Z} \psi_j(t)=1, \forall t\neq 0,
\endeq
where $\psi_j(t):=\psi_0(\lambda^{j}t)$. Hence
\beq\label{1401eee3.6}
T_A f(x)=\sum_{j\in \Z}\int_{\R}e^{iA(x)|y|^n}\psi_j(A(x)^{1/n}y)f(x- y)\frac{dy}{y}.
\endeq
The high frequency part of the kernel in \eqref{1401eee3.6}, which is 
\beq
T_A f(x)=\sum_{j\in \N}\int_{\R}e^{iA(x)|y|^n}\psi_j(A(x)^{1/n}y)f(x- y)\frac{dy}{y},
\endeq
can be bounded by
\beq
\begin{split}
& H^{*}f(x)+\sum_{j\in \N} \int_{\R}\left( e^{iA(x)|y|^n}-1\right)\psi_j(A(x)^{1/n}y)f(x- y)\frac{dy}{y}\\
& \le H^{*}f(x)+\sum_{j\in \N} \int_{2^{-\frac{j+1}{n}}A(x)^{-1/n}}^{2^{-\frac{j-1}{n}}A(x)^{-1/n}} A(x)|t|^{n-1}|f(x-t)|dt\\
& \le H^{*}f(x)+ \sum_{j\in \N} A(x) A(x)^{-\frac{n-1}{n}} 2^{-j\cdot \frac{n-1}{n}} \int_{2^{-\frac{j+1}{n}}A(x)^{-1/n}}^{2^{-\frac{j-1}{n}}A(x)^{-1/n}} |f(x-t)|dt\\
& \lesim \sum_{j\in \N} 2^{-j} Mf(x) \lesim Mf(x).
\end{split}
\endeq
Here all the constants are uniform for large $n$.\\

Concerning the low frequency part of the kernel in \eqref{1401eee3.6}, we denote
\beq
\sum_{j\ge 0}T_j f(x):=\sum_{j\ge 0} \int_{\R} f(x-t) e^{iA(x)|t|^n} \psi_j(A(x)^{1/n}t)\frac{dt}{t}.
\endeq
By applying the triangle inequality, it suffices to prove
\begin{lem}\label{1401lemma3.1}
There exists a universal constant $C>0$ and $\gamma>0$ such that for all $j\in \N$, we have
\beq
\|T_j f\|_2 \le C 2^{-\gamma j} \|f\|_2.
\endeq
\end{lem}
\noindent {\bf Proof of Lemma \ref{1401lemma3.1}:} Similar to the calculation from \eqref{1401eee2.18} to \eqref{EE3.15}, the proof of the above lemma is reduced to the following pointwise estimate of the kernel of the operator $T_j T_j^{*}$:
\begin{lem}\label{1401lemma3.2}
There exists a universal constant $C>0$ such that
\beq\label{1401eee3.11}
\begin{split}
& \left| \int_{\R}e^{i 2^{j} |\eta|^n-i 2^{j}|\xi-h \eta|^n} \frac{\psi_0(\eta)}{\eta} \frac{\psi_0(\xi- h \eta)}{\xi- h \eta} d\eta \right| \\
& \le \frac{C}{n}\chi_{[-n\cdot 2^{-j/4}, n\cdot 2^{-j/4}]}(\xi)+ C\cdot 2^{-j/4}\chi_{[-2,2]}(\xi).
\end{split}
\endeq
\end{lem}
\noindent {\bf Proof of Lemma \ref{1401lemma3.2}:} There are two cases $\xi-h\eta>0$ and $\xi-h\eta<0$. The former case remains the same as in Lemma \ref{lemma2.1}. For the latter case, denote $h_0=1/10$, then 
\beq
\text{Case } \xi-h\eta<0, 0< h<h_0
\endeq 
also remains the same. \\

{\bf Case $h\ge h_0$ and $\xi\in [-n\cdot 2^{-j/4 }, n\cdot 2^{-j/4}]$:} In this case, we bound the integrand on the left hand side of \eqref{1401eee3.11} by its absolute value to obtain
\beq
\left| \int_{\R}e^{i 2^{j} |\eta|^n-i 2^{j}|\xi-h \eta|^n} \frac{\psi_0(\eta)}{\eta} \frac{\psi_0(\xi- h \eta)}{\xi- h \eta} d\eta \right| \le 2^{1/n}-1 \le C/n,
\endeq
for some universal constant $C>0$. In this way, we obtain the first term on the right hand side of \eqref{1401eee3.11}.\\

%
%
%
%
%
%
%
%
%
%
%
%
%
%
%
%
%
%

{\bf Case $h\ge h_0$ and $|\xi| \ge n\cdot 2^{-j/4}$:} the derivative of the phase function on the left hand side of \eqref{1401eee3.11} is
\beq
n\cdot 2^j (\eta^{n-1}- h(h\eta-\xi)^{n-1}).
\endeq
The derivative might be small in this case, hence similar as before we single out a set given by
\beq
E_{h, \xi}:=\{\eta: \eta\in (1/\lambda, \lambda^2), h\eta-\xi\in (1/\lambda, \lambda^2) \text{ and } |\eta^{n-1}-h\cdot (h \eta-\xi)^{n-1}|\le n^{-1}\cdot 2^{-j/2} \},
\endeq
and what remains is to prove 
\begin{lem}\label{1401lemma3.3}
Under the above notations, we have $|E_{h,\xi}|\lesim 2^{-j/4}.$
\end{lem}

\noindent {\bf Proof of Lemma \ref{1401lemma3.3}:} By the fundamental theorem, we obtain
\beq
|\eta^{n-1}-h\cdot (h \eta-\xi)^{n-1}| \sim n |\eta- h^{\frac{n}{n-1}}\eta+ h^{\frac{1}{n-1}}\xi|.
\endeq
Hence any point $\eta\in E_{h,\eta}$ satisfies
\beq\label{29novEE3.16}
n |\eta- h^{\frac{n}{n-1}}\eta+ h^{\frac{1}{n-1}}\xi| \lesim n^{-1}\cdot 2^{-j/2} .
\endeq
This implies 
\beq
|E_{h,\xi}|\lesim n^{-2}\cdot 2^{-j/2}(1-h^{\frac{n}{n-1}})^{-1}.
\endeq
In order for the set $E_{h,\xi}$ not to be empty, we need 
\beq
|1-h^{\frac{n}{n-1}}| \ge 2^{-j/4},
\endeq
as otherwise the inequality \eqref{29novEE3.16} will not hold true. Hence
\beq
|E_{h,\xi}|\lesim n^{-2} 2^{-j/4}.
\endeq
So far we have finished the proof of Lemma \ref{1401lemma3.3}, thus the proof of the uniform estimate in Theorem \ref{1401theorem1.3}.

\section{Proof of Theorem \ref{1401theorem1.6}}\label{1602section4}

In this section, we present the proofs of the uniform estimates \eqref{1401eee1.17} and \eqref{1401eee1.18}. This time, unlike the situation for Theorem \ref{theorem1.1} and Theorem \ref{1401theorem1.3}, the arguments for these two proofs are no longer similar, hence we present them in the following two subsections separately.

\subsection{Proof of the estimate \eqref{1401eee1.17}}

In this subsection we will prove the first part of Theorem \ref{1401theorem1.6}: There exists a universal constant $C$ such that for all $\lambda\in \R$ and all $|\epsilon|\le 1/2$:
\beq
|\int_{\R} e^{i|t|^{\epsilon}} e^{-i\lambda t}\frac{dt}{t}| \le C.
\endeq
In the following, again we will only write down the proof for positive $\epsilon$. The proof for negative $\epsilon$ is similar.\\

By the change of variable $t\to -t$, it is clear that we only need to look at the case $\lambda>0$. After another change of variable
\beq
\lambda t\to t,
\endeq 
it suffices to prove the uniform bound
\beq\label{812EE2.3}
|\int_{\R} e^{i \lambda |t|^{\epsilon}-it}\frac{dt}{t}|\le C.
\endeq

Notice that the function $\frac{e^{i|t|^{\epsilon}}}{t}$ is an odd function, hence the integration of this function over $\R$ is zero. However, there is still another part $e^{-it}$ in the phase function, which makes the integrand no longer odd. The idea is to approximate $e^{it}$ by constant 1 when $t$ is small. We split the integration in \eqref{812EE2.3} into the following two parts:
\beq
\int_{\R} e^{i \lambda |t|^{\epsilon}-it}\frac{dt}{t}=\int_{-1}^1 e^{i \lambda |t|^{\epsilon}-it}\frac{dt}{t}+\int_{\R\setminus [-1, 1]} e^{i \lambda |t|^{\epsilon}-it}\frac{dt}{t}.
\endeq
We denote the first term by $I$, and the second term by $II$. For the first term:
\beq
\begin{split}
I&= \int_{-1}^0 e^{i \lambda |t|^{\epsilon}-it}\frac{dt}{t}+\int_{0}^1 e^{i \lambda |t|^{\epsilon}-it}\frac{dt}{t}\\
	&= \int_0^1 (e^{i\lambda |t|^{\epsilon}-it}-e^{i\lambda |t|^{\epsilon}+it})\frac{dt}{t}.
\end{split}
\endeq
Hence
\beq
|I| \lesim \int_0^1 |e^{-it}-e^{it}|\frac{dt}{|t|} \lesim 1.
\endeq
For the second term $II$, we first write it as
\beq\label{812EE2.7}
II= \int_1^{\infty} e^{i \lambda t^{\epsilon}+it}\frac{dt}{t}+\int_1^{\infty} e^{i \lambda t^{\epsilon}-it}\frac{dt}{t}.
\endeq
For the former term, we see that the phase function $\lambda t^{\epsilon}+t$ does not have any critical point on the interval $[1, \infty]$, which suggests that this term can simply be bounded by doing integration by part:
\beq
\left| \int_1^{\infty} e^{i \lambda t^{\epsilon}+it}\frac{dt}{t} \right| = \left| \int_1^{\infty} \frac{1}{\lambda \epsilon t^{\epsilon}+t} d(e^{i\lambda t^{\epsilon}+it}) \right| \lesim 1.
\endeq
For the latter term in \eqref{812EE2.7}, whether the phase function $\lambda t^{\epsilon}-t$ has critical point or not depends on the choice of the parameters $\lambda$ and $\epsilon$.

\subsubsection{The case $\lambda \cdot \epsilon \le 1/10$}

In this case, it is not difficult to see that the phase function $\lambda t^{\epsilon}-t$ in the latter term of \eqref{812EE2.7} has no critical point. Hence it suffices to do an integration by part:
\beq
\left|\int_1^{\infty} e^{i\lambda t^{\epsilon}-it}\frac{dt}{t}\right| = \left| \int_1^{\infty} \frac{1}{\lambda \epsilon t^{\epsilon}-t}d(e^{i\lambda t^{\epsilon}-t}) \right| \lesim 1.
\endeq

\subsubsection{The case $\lambda \cdot \epsilon \ge 1/10$}\label{subsection2.2}

The term we need to bound is
\beq\label{812EE2.10}
\int_1^{\infty} e^{i\lambda t^{\epsilon}-it}\frac{dt}{t}.
\endeq
The phase function $\lambda t^{\epsilon}-t$ has a critical point at 
\beq
t_0:=\left( \lambda \epsilon \right)^{\frac{1}{1-\epsilon}}.
\endeq
Hence we split the integration in \eqref{812EE2.10} into two parts accordingly:
\beq\label{812EE2.12}
\int_1^{\infty} e^{i\lambda t^{\epsilon}-it}\frac{dt}{t}= \int_1^{t_0} e^{i\lambda t^{\epsilon}-it}\frac{dt}{t}+ \int_{t_0}^{\infty} e^{i\lambda t^{\epsilon}-it}\frac{dt}{t}.
\endeq
We denote the former term in the last expression by $III$ and the latter term by $IV$. \\

For the term $III$, the function $e^{i\lambda t^{\epsilon}}$ has higher oscillation than the function $e^{it}$, which suggests the following integration by part:
\beq
\begin{split}
III&= \left| \int_1^{t_0} \frac{1}{\lambda \epsilon t^{\epsilon}} e^{it}d(e^{i\lambda t^{\epsilon}}) \right| \lesim 1+ \left| \int_1^{t_0} e^{i\lambda t^{\epsilon}} \left( \frac{e^{it}}{\lambda \epsilon t^{\epsilon}}\right)' dt\right|\\
	& \lesim 1+ \int_1^{t_0} \left(\frac{1}{\lambda \epsilon t^{\epsilon}}+ \frac{1}{\lambda t^{\epsilon+1}}\right)dt \lesim 1+ \int_1^{t_0} \frac{1}{\lambda \epsilon t^{\epsilon}}dt  \lesim 1.
\end{split}
\endeq
For the term $IV$, the roles of the two functions $e^{i\lambda t^{\epsilon}}$ and $e^{it}$ are reversed:
\beq
\begin{split}
IV & = \left| \int_{t_0}^{\infty}\frac{e^{i\lambda t^{\epsilon}}}{t} d(e^{it})  \right|\lesim 1+ \left| \int_{t_0}^{\infty} e^{it}\left(\frac{e^{i\lambda t^{\epsilon}}}{t} \right)' dt\right|\\
		&\lesim 1 + \int_{t_0} ^{\infty} \left| \frac{i\lambda \epsilon t^{\epsilon} e^{i\lambda t^{\epsilon}}-e^{i\lambda t^{\epsilon}}}{t^2} \right|dt \lesim 1+ \lambda \epsilon \int_{t_0} ^{\infty} \frac{1}{t^{2-\epsilon}}dt \lesim 1.
\end{split}
\endeq
So far we have finished the proof for the case $\lambda \epsilon\ge 1/10$, thus the first part of Theorem \ref{1401theorem1.6}.

\subsection{Proof of the estimate \eqref{1401eee1.18}}

In this subsection, we will show that there exists a universal constant $C>0$ such that for all $1/2<\epsilon<3/2$ and all $\lambda\in \R$, we have 
\beq
\left| \int_{\R} e^{i \lambda\cdot  \text{sgn}(t) |t|^{\epsilon}-it}\frac{dt}{t}\right| \le C.
\endeq
First notice that by doing the change of variable $\lambda \cdot \text{sgn}(t)\cdot |t|^{\epsilon}\to t$, it suffices to consider the case $\epsilon>1$. We further simplify the above estimate by using some trivial cancellation:
\beq\label{812EE3.2}
\begin{split}
& \int_{\R} e^{i \lambda\cdot  \text{sgn}(t) t^{\epsilon}-it}\frac{dt}{t}\\
&= \int_0^{\infty} e^{i \lambda t^{\epsilon}-it}\frac{dt}{t}+\int_{-\infty}^0 e^{-i\lambda |t|^{\epsilon}-it}\frac{dt}{t}=\int_0^{\infty} \sin (\lambda t^{\epsilon}-t)\frac{dt}{t}.
\end{split}
\endeq

For the case $\lambda<0$, we see easily that there is no critical point of the phase function $\lambda t^{\epsilon}-t$. Hence this case is supposed to be easier: Take $t_0$ such that
\beq
-\lambda t_0^{\epsilon}+t_0=1.
\endeq
We split the integration in the last term of \eqref{812EE3.2} into two parts:
\beq\label{812ee3.4}
\left(\int_0^{t_0}+\int_{t_0}^{\infty}\right) \sin(-i\lambda t^{\epsilon}+it)\frac{dt}{t}.
\endeq
To bound the former part of the last expression, the idea is to use the simple inequality that $|\sin t | \le |t|$ when $t$ is small:
\beq
\left| \int_0^{t_0}\sin (-i\lambda t^{\epsilon}+it) \frac{dt}{t}\right| \lesim \int_0^{t_0} \left(-\lambda t^{\epsilon-1}+1\right) dt \le -\frac{\lambda}{\epsilon}t_0^{\epsilon}+t_0 \lesim \frac{1}{\epsilon}.
\endeq
For the latter part of \eqref{812ee3.4}, we will do an integration by part to explore the high oscillation from the term $\sin (-i\lambda t^{\epsilon}+it)$:
\beq
\begin{split}
& \left|\int_{t_0}^{\infty} \sin (-i\lambda t^{\epsilon}+it)\frac{dt}{t}\right|\\
& \lesim  \left|\int_{t_0}^{\infty} \frac{1}{-\lambda \epsilon t^{\epsilon}+t}d(\sin (-i\lambda t^{\epsilon}+t))\right| \lesim \frac{1}{-\lambda \epsilon t_0^{\epsilon}+t_0} \lesim 1.
\end{split}
\endeq
So far we have finished the proof of the case $\lambda <0$. In the following, we will focus on the case $\lambda>0$. Moreover, we will write 
\beq
\epsilon=1+\frac{1}{n}
\endeq
from time to time to indicate that $n$ is some large number.

\subsubsection{The case $\lambda>n^{-1/n}$}

In this case the minimum of the phase function $\lambda t^{\epsilon}-t$ is 
\beq
-\left(\frac{1}{\lambda}\right)^n \frac{1}{n}\cdot \left(\frac{n}{n+1}\right)^{n+1}\ge -1.
\endeq
We denote by $t_0$ such that
\beq
\lambda t_0^{\epsilon}-t_0 =0,
\endeq
and by $t_2$ such that
\beq
\lambda t_2^{\epsilon}-t_2 =1.
\endeq
Observe that when $t\le t_2$, the absolute value of the phase function is small, i.e.
\beq
|\lambda t^{\epsilon}-t|\le 1.
\endeq
This suggests the following splitting of the term \eqref{812EE3.2} that we need to bound:
\beq
\left(\int_0^{t_2}+\int_{t_2}^{\infty}\right) \sin (\lambda t^{\epsilon}-t)\frac{dt}{t}.
\endeq
The former part will be denoted by $V$, and the latter part $VI$. For the term $VI$, we simply do an integration by part:
\beq
|VI| \lesim \left|\int_{t_2}^{\infty}\frac{1}{\lambda \epsilon t^{\epsilon}-t}d (\sin (\lambda t^{\epsilon}-t))\right| \lesim 1.
\endeq
For the term $V$, we use the simple inequality that $|\sin t|\le |t|$ for small $t$:
\beq\label{812ee3.28}
\begin{split}
|V| &\lesim \int_0^{t_2} |\lambda t^{\epsilon-1}-1|dt\\
	& \lesim \int_0^{t_0} (-\lambda t^{\epsilon-1}+)dt+ \int_{t_0}^{t_2} (\lambda t^{\epsilon-1}-1)dt\\
	& \lesim \left(t_0-\frac{\lambda}{\epsilon}t_0^{\epsilon}\right) + \left( \frac{\lambda}{\epsilon}t_2^{\epsilon}-t_2-\frac{\lambda}{\epsilon}t_0^{\epsilon}+t_0 \right).
\end{split}\endeq
For the latter part of the last expression in \eqref{812ee3.28}, by the definition of $t_0$ and $t_2$, we obtain that 
\beq\label{812ee3.29}
\frac{\lambda}{\epsilon}t_2^{\epsilon}-t_2-\frac{\lambda}{\epsilon}t_0^{\epsilon}+t_0 =\frac{1}{\epsilon}-(1-\frac{1}{\epsilon})(t_2-t_0).
\endeq
We know \eqref{812ee3.29} must positive as the integrand is positive, hence \eqref{812ee3.29} can be bounded by $\frac{1}{\epsilon}\le 1$. For the former part of the last expression in \eqref{812ee3.28}, by the definition of $t_0$, we obtain
\beq
t_0-\frac{\lambda}{\epsilon}t_0^{\epsilon}\lesim  \frac{t_0}{n}\le \left(\frac{1}{\lambda}\right)^n \frac{1}{n}\le 1.
\endeq
So far we have finished the proof of the case $1\le \epsilon<3/2, \lambda>n^{-1/n}$.

\subsubsection{The case $\lambda \le n^{-1/n}$}

We denote by $t_2$ the smaller one of the two positive numbers such that
\beq
\lambda t_2^{\epsilon}-t_2=-2.
\endeq
For $t\le t_2$, again we observe the fact that the phase function is small, which suggests to write
\beq
\eqref{812EE3.2}=\left( \int_0^{t_2}+ \int_{t_2}^{\infty} \right) \sin (\lambda t^{\epsilon}-t)\frac{dt}{t}.
\endeq
We denote the first term of the last expression by $VII$, and the second by $VIII$. To estimate the term $VII$, we do the following routine calculation: 
\beq
|VII|\lesim \int_0^{t_2}\left( 1-\lambda t^{\epsilon-1} \right)dt \le t_2-\frac{\lambda}{\epsilon}t_2^{\epsilon}=\frac{t_2}{n}+\frac{1}{\epsilon}.
\endeq
To finish the estimate of the term $VII$, we need the following 
\begin{lem}\label{1401lemma4.1}
Under the above notations, we have $t_2\le n$.
\end{lem}
\noindent {\bf Proof of Lemma \ref{1401lemma4.1}:} If we fix one $\epsilon$, then $t_2$ can be viewed as a function of $\lambda$. Moreover, it is easy to see that this function is monotone increasing. Hence we only need to prove that $t_2(n^{-1/n})\le n$, which is trivial. $\Box$\\

For the term $VIII$,   we would like to do the following integration by part:
\beq
VIII=\int_{t_2}^{\infty} \frac{1}{\lambda \epsilon t^{\epsilon}-t}d(\cos (\lambda t^{\epsilon}-t)).
\endeq
However, notice that the denominator $\lambda \epsilon t^{\epsilon}-t$ is not aways small, or in another word, the phase function $\lambda t^{\epsilon}-t$ does not always oscillate fast on the interval $(t_2, \infty)$. Hence we need to do a finer decomposition for the interval $(t_2, \infty)$.\\

Denote by $t_3$ such that
\beq
\lambda \epsilon t_3^{\epsilon}-t_3=-1,
\endeq
and by $t_4$ that 
\beq
\lambda \epsilon t_4^{\epsilon}-t_4=1.
\endeq
We split the integration in $VIII$ into the following:
\beq
VIII= \left(\int_{t_2}^{t_3}+\int_{t_3}^{t_4}+\int_{t_4}^{\infty}\right) \sin (\lambda t^{\epsilon}-t)\frac{dt}{t}.
\endeq
Simply by integration by part, we obtain
\beq
|VIII| \lesim \frac{1}{|\lambda \epsilon t_2^{\epsilon}-t_2|} + \frac{t_4}{t_3}.
\endeq
Notice that
\beq
\lambda \epsilon t_2^{\epsilon}-t_2=\frac{\lambda t_2^{\epsilon}}{n}-2=\frac{t_2-2}{n}-2\le -1.
\endeq
Hence what is left is to prove the following 
\begin{lem}\label{812lemma3.2}
Under the above notations, we have that $t_4\le 10^3 t_3$.
\end{lem}
\noindent {\bf Proof of Lemma \ref{812lemma3.2}:} by definition, we have
\beq
\lambda \epsilon t_3^{\epsilon}-t_3=-1.
\endeq
By monotonicity of the function $\lambda \epsilon t^{\epsilon}-t$ for $t\ge t_3$, to prove that $t_4\le 10^3 t_3$, it suffices to prove that 
\beq
\lambda \epsilon (10^3 t_3)^{\epsilon}-10^3 t_3\ge 1.
\endeq
We substitute the definition of $t_3$ into the last expression to obtain
\beq\label{812ee3.42}
\lambda \epsilon (10^3 t_3)^{\epsilon}-10^3 t_3=\lambda\epsilon(10^{3(1+\frac{1}{n})}-10^3)t_3^{\epsilon}-10^3 \ge 10^3\ln 10^3 \frac{t_3-1}{n}-10.
\endeq
To show that the last expression in \eqref{812ee3.42} is greater than 1, it suffices to prove
\begin{lem}\label{812lemma3.3}
Under the above notations, we have that $t_3\ge n/10.$
\end{lem}
\noindent {\bf Proof of Lemma \ref{812lemma3.3}:} The minimum of the function $\lambda \epsilon t^{\epsilon}-t$ is attained at 
\beq 
t=\left(\frac{1}{\lambda \epsilon^2}\right)^n\ge n \left( \frac{1}{\epsilon}\right)^{2n}=n \left( \frac{n}{n+1}\right)^{2n}\ge n/10.
\endeq
By definition, $t_3$ lies on the right hand side of the critical point of the function $\lambda \epsilon t^{\epsilon}-t$. Hence $t_3\ge n/10.$ So far we have finished the proof of Lemma \ref{812lemma3.3}, hence Lemma \ref{812lemma3.2}, thus Theorem \ref{1401theorem1.6}. $\Box$

%
%
%
%
%
%
%
%
%
%
%
%
%
%
%
%
%
%

Shaoming Guo, Institute of Mathematics, University of Bonn\\
\indent Address: Endenicher Allee 60, 53115, Bonn\\
\indent Email: shaoming@math.uni-bonn.de

\end{document}